\documentclass[12pt]{article}
\usepackage{latexsym, amssymb, amsmath, amscd, amsfonts, epsfig, graphicx, colordvi,verbatim,ifpdf,tikz,extarrows}
\usepackage{amsfonts, amsmath, amssymb}
\usepackage{amssymb,amsfonts,amsmath,latexsym,epsfig,cite, psfrag,eepic,color}
\usepackage{amscd,graphics}
\usepackage{latexsym, amssymb,  amsmath,amscd, amsfonts, epsfig, graphicx, colordvi,amsthm}
\usepackage{ytableau}

\usepackage[colorlinks,linkcolor=blue,anchorcolor=blue,citecolor=blue]{hyperref}

\usepackage{graphicx}
\usepackage{color}
\usepackage{ifpdf}
\usepackage{fancybox}

\usepackage{float}

\newtheorem{thm}{Theorem}[section]

\newtheorem{conj}[thm]{Conjecture}

\newtheorem{lem}[thm]{Lemma}

\def\pf{\noindent{\it Proof.} }
\def\qed{\nopagebreak\hfill{\rule{4pt}{7pt}}
\medbreak}

\hypersetup{colorlinks = true}

\numberwithin{equation}{section}

\allowdisplaybreaks

\def\qed{\nopagebreak\hfill{\rule{4pt}{7pt}}
\medbreak}

\newlength{\boxedparwidth}
  {\begin{center} \begin{tabular}{|@{\hspace{.315in}}c@{\hspace{.15in}}|}
                  \hline \\ \begin{minipage}[t]{\boxedparwidth}
                  \setlength{\parindent}{.25in}}%
  {\end{minipage} \\ \\ \hline \end{tabular} \end{center}}

\begin{document}
\begin{center}

 {\Large \bf Proof of Singh and Barman's conjecture on hook length biases}
\end{center}

\begin{center}
  {Hongshu Lin$^1$}, and {Wenston J.T. Zang$^2$} \vskip 2mm

 $^{1,2}$School of Mathematics and Statistics, Northwestern Polytechnical University, Xi'an 710072, P.R. China\\[6pt]
	$^{1,2}$ MOE Key Laboratory for Complexity Science in Aerospace, Northwestern Polytechnical University, Xi'an 710072, P.R. China\\[6pt]
	$^{1,2}$ Xi'an-Budapest Joint Research Center for Combinatorics, Northwestern Polytechnical University, Xi'an 710072, P.R. China\\[6pt]

   \vskip 2mm

$^1$linhongshu@mail.nwpu.edu.cn, 
    $^2$zang@nwpu.edu.cn
\end{center}

\vskip 6mm \noindent {\bf Abstract.} Let $b_{t,i}(n)$ denote the total number of $i$-hooks in $t$-regular partitions of $n$. Singh and Barman conjectured that $b_{t+1,2}(n) \geq b_{t,2}(n)$ holds for all $t\ge 3$ and $n\ge 0$. This conjecture was known to hold for $t=3$ due to work of Barman Mahanta and Singh. In this paper, we prove this conjecture.



\section{Introduction}


This paper investigates the $2$-hook length biases in $t$-regular partitions and $t+1$-regular partitions. Recall that a partition of a positive integer $n$ is a finite sequence of non-increasing positive integers $( \lambda_{1},\lambda_{2},\dots,\lambda_{\ell} ) $ satisfying the condition that the sum of its terms $\lambda_{1} + \lambda_{2}+ \dots + \lambda_{\ell} = n $. Given $t \geq 2$, a $t$-regular partition is a partition where none of its parts is divisible by $t$, and a $t$-distinct partition is a partition where each number appears at most $t-1$ times. We denote by $b_t(n)$ the number of $t$-regular partitions of $n$ and denote by $d_t(n)$ the number of $t$-distinct partitions of $n$. For example, given $n=6$ and $t=2$, there are four 2-regular partitions of $6$, namely
\[\begin{array}{cccc}
(5,1), & (3,3), & (3,1,1,1), & (1,1,1,1,1,1).
\end{array}\]
And there are four 2-distinct partitions of $6$, namely
\[\begin{array}{cccc}
(6), & (5,1), & (4,2), & (3,2,1).
\end{array}\]
Thus $b_2(6)=d_2(6)=4$.

\begin{figure}[h]
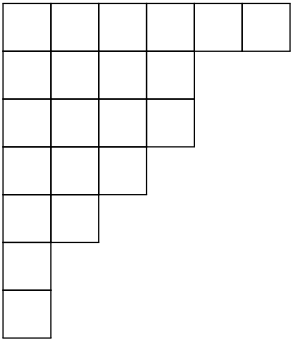

		\begin{center}
		\begin{ytableau}[]
        \ & \ & \ & \ & \ & \ \\
        \ & \ & \ & \ \\
        \ & \ & \ & \ \\
        \ & \ & \  \\
        \ & \ \\
        \ \\
        \
        \end{ytableau}
		\caption{The Young diagram of partition $(6,4,4,3,2,1,1)$}
		\label{F1}
		\end{center}
	\end{figure}

\begin{figure}[h]
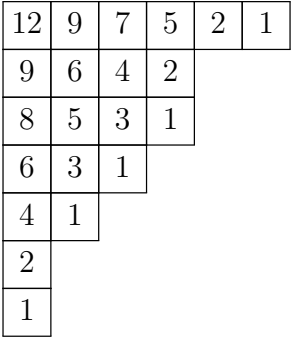

		\begin{center}
		\begin{ytableau}[]
        12 & 9 & 7 & 5 & 2 & 1 \\
        9 & 6 & 4 & 2 \\
        8 & 5 & 3 & 1 \\
        6 & 3 & 1 \\
        4 & 1 \\
        2 \\
        1
        \end{ytableau}
		\caption{The hook lengths of $(6,4,4,3,2,1,1)$}
		\label{F2}
		\end{center}
	\end{figure}

Hook length is one of the basic definition in the theory of integer partitions. Recall that the Young diagram corresponding to a partition $( \lambda_{1}, \lambda_{2}, \dots, \lambda_{\ell})$ is a left-justified array of boxes, where the $i$-th row (counting from the top) contains exactly $\lambda_{i}$ boxes. For any given box in a Young diagram, its hook length is defined as the sum of three components: the number of boxes that lie directly to the right of the box in question, the number of boxes that lie directly below it, and $1$ (this $1$ accounts for the box itself). For example, the Young diagram of partition $(6,4,4,3,2,1,1)$ is shown in Figure \ref{F1}. The hook length of a box in the Young diagram is the number of the boxes directly to its right or directly below it and including itself exactly once. For example, Figure \ref{F2} shows the hook length of each box in the Young diagram of partition $(6,4,4,3,2,1,1)$.


The hook length plays a crucial role in the representation theory of the symmetric group $S_n$ and the general linear group $GL_n(\mathbb{C})$. It is also closely connected to the theory of symmetric functions, particularly through the enumeration of standard Young tableaux, as formalized by the well-known hook-length formula (for further details, see \cite{Frame-Rhobinson-Thrall-1954, Young-1927}). Furthermore, the Nekrasov-Okounkov formula (refer to \cite{Han-Nek-2010, Nek-Oko}) establishes a fundamental link between hook length and the theory of modular forms. Beyond these core connections, extensive research has been devoted to studying the properties of hook length, including its asymptotic behavior, combinatorial characteristics, and arithmetic properties. Notable among these studies are investigations into $t$-core partitions and $t$-hook partitions; relevant references include \cite{Garvann-Kim, James-Kerber, Littlewood-Modular}.

 Let $b_{t,k}(n)$ denote the number of hooks of length $k$ across all $t$-regular partitions of $n$ and $d_{t,k}(n)$.   Ballantine, Burson, Craig, Folsom and Wen \cite{Ballantine-2023} investigate the bias in the number of hooks of fixed lengths between odd partitions and distinct partitions. Craig, Dawsey and Han \cite{Craig-2023} found the generating function and asymptotics result on the number of $i$-hooks in odd parts and distinct parts. Kim \cite{Kim} established generating function and asymptotic results on the total number of hooks of lengths 1, 2, 3 in $\ell$-regular and $\ell$-distinct partitions.  For the fixed integer $t=2$, Singh and Barman \cite{Singh-Barman-2024} found two key inequalities: $b_{2,2}(n) \geq b_{2,1}(n)$ for all $n > 4$, and $b_{2,2}(n) \geq b_{2,3}(n)$ for all $n \geq 0$. In addition, they proposed the conjecture that $b_{3,2}(n) \geq b_{3,1}(n)$ holds for all $n \geq 28$; this conjecture was recently verified by He and Liu \cite{He-liu-2025}, Qu and Zang \cite{qu-zang-2025} independently.

This paper focuses on investigating the inequality between $b_{t+1,2}(n) $ and $b_{t,2}(n)$.  Singh and Barman \cite{Singh-Barman-2024} established that the inequality $b_{t+1,1}(n) \geq b_{t,1}(n)$ holds for all integers $t \geq 2$ and $n \geq 0$. In a subsequent work \cite{Singh-Barman-2025}, Singh and Barman initiated the study of $2$-hook biases between $t$-regular partitions and $t+1$-regular partitions as follows.
\begin{thm}[\cite{Singh-Barman-2025}]\label{thm-sb-2025}
    For all integers $n > 3$, we have $b_{3,2}(n) \geq b_{2,2}(n)$.
\end{thm}
Singh and Barman \cite{Singh-Barman-2025} also proposed the following conjecture concerning the relation between $b_{t+1,2}(n)$ and $b_{t,2}(n)$:
\begin{conj}[\cite{Singh-Barman-2025}]\label{conj1}
     Let $t \geq 3$ be an integer. We have $b_{t+1,2}(n) \geq b_{t,2}(n)$ for all $n \geq 0$.
\end{conj}
 Barman, Mahanta and   Singh \cite{Barman-2025} verified the validity of this conjecture for the specific case $t = 3$ as given below. 

\begin{thm}[\cite{Barman-2025}]
    Conjecture \ref{conj1} is true for $t = 3$.
\end{thm}
Kim \cite{Kim} confirmed this conjecture for sufficiently large $n$.
\begin{thm}[\cite{Kim}]
    Let $t\geq 2$ be an integer. For each $\ell=1,2,3$, there exist some positive integers $N_t$ such that for all $n>N_t$,
    \begin{equation}
        b_{t+1,\ell}(n)\geq b_{t,\ell}(n).
    \end{equation}
\end{thm}

The main result of this paper is the following theorem.
\begin{thm}\label{thm-com-int-umn}
   Conjecture \ref{conj1} is true.
\end{thm}

Combining Theorem \ref{thm-sb-2025}, we obtain the following generalized result.

\begin{thm}
    For integers $t\geq 2$ and $n\geq0$ except for $(t,n)=(2,3)$, we have
    \begin{equation}
        b_{t+1,2}(n) \geq b_{t,2}(n).
    \end{equation}
\end{thm}

This paper is organized as follows. Section 2 is devoted to proving Theorem \ref{thm-com-int-umn} for odd $t$. Specifically, we begin by decomposing the generating function for $b_{t+1,2}(n) - b_{t,2}(n)$ into six distinct parts. By constructing three injections between these parts, we establish the non-negativity of the generating function in the case of odd $t$. Section 3 addresses the case of even $t$, where a similar proof strategy is employed for $t \geq 6$. The situation for $t = 4$, however, requires a minor modification in the third injection. With this adjustment, the proof of Theorem \ref{thm-com-int-umn} is completed.

\section{Proof of Theorem \ref{thm-com-int-umn} for odd number $t$}\label{2}
This section proves Theorem \ref{thm-com-int-umn} for odd $t$. We begin by decomposing the generating function for $b_{t+1,2}(n) - b_{t,2}(n)$ into six components and providing a combinatorial interpretation for each. The proof is then completed by constructing three key injections.

We begin by introducing some notation on integer partitions. Throughout this paper, for a partition $\lambda$ of $n$, let $f_\lambda(k)$ denote the number of times $k$ appears in $\lambda$. With this notation, we write $\lambda$ as $(n^{f_\lambda(n)},(n-1)^{f_\lambda(n-1)},\ldots,1^{f_\lambda(1)})$, where when $f_\lambda(i)=0$, we may omit the term $i^0$, and we write $i^1$ as $i$ for simplicity. For example, $\lambda=(5,3,3,1,1,1)$. We have that $f_\lambda(5)=1$, $f_\lambda(3)=2$, $f_\lambda(1)=3$ and $f_\lambda(i)=0$ for all $i\ne 1,3,5$. We can write $\lambda=(5,3^2,1^3)$ instead of $(5,3,3,1,1,1)$.
Let $\lambda,\mu$ be two partitions. We use $\lambda\cup\mu$ to denote the partition
\[(m^{f_\lambda(m)+f_\mu(m)},(m-1)^{f_\lambda(m-1)+f_\mu(m-1)},\ldots,1^{f_\lambda(1)+f_\mu(1)}),\]where $m=\max\{\lambda_1,\mu_1\}$.
Moreover, if $f_\lambda(k)\ge f_\mu(k)$ for all $k$, then we denote the partition
\[(n^{f_\lambda(n)-f_\mu(n)},(n-1)^{f_\lambda(n-1)-f_\mu(n-1)},\ldots,1^{f_\lambda(1)-f_\mu(1)})\]
by $\lambda\setminus\mu$. For example, $\lambda=(6,5^2,2^4,1^5)$ and $\mu=(5,2^3,1^2)$, we have $\lambda\cup\mu=(6,5^3,2^7,1^7)$ and $\lambda\setminus\mu=(6,5,2,1^3)$.

Corteel and Lovejoy \cite{Corteel-Lovejoy-2004} first introduced the concept of an overpartition, which is defined as a partition in which the first occurrence of any given part may be marked. For example, there are $14$ overpartitions of $4$, namely
\[\begin{array}{cccccccccc}
     (4),&(\overline{4}),&(3,1),&(\overline{3},1),&(3,\overline{1}),&(\overline{3},\overline{1}),&(2,2),  \\
     (\overline{2},2),&(2,1,1),&(\overline{2},1,1),&(2,\overline{1},1),&(\overline{2},\overline{1},1),&(1,1,1,1),&(\overline{1},1,1,1).
\end{array}\]

In this paper, in sake of convenience, let OPO-overpartition denote the overpartition with exactly one part overlined. In the above example, there are $7$ OPO-overpartition of $4$, which is $(\overline{4}), (\overline{3},1), (3,\overline{1}), (\overline{2},2), (\overline{2},1,1), (2,\overline{1},1), (\overline{1},1,1,1)$. Given an OPO-overpartition $\lambda$, we use $o(\lambda)$ to denote the unique overlined part. For example, for the OPO-overpartition $\lambda=(\overline{3},1)$, $o(\lambda)=3$. Moreover, we further define a $t$-regular  OPO-overpartition as a OPO-overpartition with each part not a multiple of $t$. It is clear that there are three $2$-regular OPO-overpartition of $4$, namely $(\overline{3},1), (3,\overline{1}), (\overline{1},1,1,1)$.

We first recall the generating function of $b_{t,2}(n)$ given by Kim\cite{Kim}.

\begin{thm}[\cite{Kim}]\label{bt2}
    For $t\geq 2$, we have 
    \begin{equation}\label{equ-bt2}
\sum_{n=0}^\infty b_{t,2}(n)q^n =\frac{(q^t;q^t)_\infty}{(q;q)_\infty}\left(\frac{2q^2}{1-q^2}-\frac{q^t}{1-q^t}+\frac{q^{2t-1}-q^{2t}+q^{2t+1}}{1-q^{2t}}\right).
    \end{equation}
\end{thm}

In order to transform the generating function of $b_{t+1,2}(n)-b_{t,2}(n)$, we give several definitions as follows. Set 
\begin{align*}
    a_{t}&=\frac{(q^{t+1};q^{t+1})_\infty}{(q;q)_\infty}\frac{\sum_{n=1}^{t}{q^{2n}}-q^{t+1}}{1-q^{2t+2}},\\
    b_t&=\frac{(q^{t};q^{t})_\infty}{(q;q)_\infty}\frac{\sum_{n=1}^{t-1}{q^{2n}}}{1-q^{2t}},\\
    c_t&=\frac{(q^{t};q^{t})_\infty}{(q;q)_\infty}\frac{q^t}{1-q^{2t}},\\
    d_{t}&=\frac{(q^{t+1};q^{t+1})_\infty}{(q;q)_\infty}\frac{\sum_{n=1}^{t}{q^{2n}}}{1-q^{2t+2}},\\
    e_{t}&=\frac{(1-q^{t+1})(q^{3t+3};q^{t+1})_\infty}{(q;q)_\infty}(q^{2t+1}+q^{2t+3}),\\
    f_t&=\frac{(1-q^{t})(q^{3t};q^{t})_\infty}{(q;q)_\infty}(q^{2t-1}+q^{2t+1}).
\end{align*}

\begin{lem}
    For given odd integer $t$, we have that 
    \begin{equation}\label{bt+1-bt}
     \sum_{n=1}^{\infty}b_{t+1,2}(n)q^n-\sum_{n=1}^{\infty}b_{t,2}(n)q^n=a_{t}-2b_t+c_t+d_{t}+e_{t}-f_{t}.
    \end{equation}

\end{lem}
\pf By Theorem \ref{bt2}, we see that 
\begin{align}
    &\sum_{n=1}^{\infty}b_{t+1,2}(n)q^n-\sum_{n=1}^{\infty}b_{t,2}(n)q^n\nonumber\\
    =&\frac{(q^{t+1};q^{t+1})_\infty}{(q;q)_\infty}\left(\frac{2q^2}{1-q^2}-\frac{q^{t+1}}{1-q^{t+1}}+\frac{q^{2t+1}-q^{2t+2}+q^{2t+3}}{1-q^{2t+2}}\right)\nonumber\\
    &-\frac{(q^t;q^t)_\infty}{(q;q)_\infty}\left(\frac{2q^2}{1-q^2}-\frac{q^t}{1-q^t}+\frac{q^{2t-1}-q^{2t}+q^{2t+1}}{1-q^{2t}}\right)\nonumber\\
    =&\frac{(q^{t+1};q^{t+1})_\infty}{(q;q)_\infty}\left(\frac{2q^2}{1-q^2}-\frac{q^{t+1}}{1-q^{t+1}}-\frac{q^{2t+2}}{1-q^{2t+2}}\right)-\frac{(q^t;q^t)_\infty}{(q;q)_\infty}\left(\frac{2q^2}{1-q^2}-\frac{q^t}{1-q^t}-\frac{q^{2t}}{1-q^{2t}}\right)\nonumber\\
    &+\frac{(q^{t+1};q^{t+1})_\infty}{(q;q)_\infty}\left(\frac{q^{2t+1}+q^{2t+3}}{1-q^{2t+2}}\right)-\frac{(q^t;q^t)_\infty}{(q;q)_\infty}\left(\frac{q^{2t-1}+q^{2t+1}}{1-q^{2t}}\right)\nonumber\\
    =&\frac{(q^{t+1};q^{t+1})_\infty}{(q;q)_\infty}\left(\frac{q^2}{1-q^2}-\frac{q^{t+1}}{1-q^{t+1}}\right)+\frac{(q^{t+1};q^{t+1})_\infty}{(q;q)_\infty}\left(\frac{q^2}{1-q^2}-\frac{q^{2t+2}}{1-q^{2t+2}}\right)\nonumber\\
    &-\frac{(q^t;q^t)_\infty}{(q;q)_\infty}\left(\frac{q^2}{1-q^2}-\frac{q^t}{1-q^t}\right)-\frac{(q^t;q^t)_\infty}{(q;q)_\infty}\left(\frac{q^2}{1-q^2}-\frac{q^{2t}}{1-q^{2t}}\right)\nonumber\\
    &+\frac{(1-q^{t+1})(q^{3t+3};q^{t+1})_\infty}{(q;q)_\infty}(q^{2t+1}+q^{2t+3})-\frac{(1-q^{t})(q^{3t};q^{t})_\infty}{(q;q)_\infty}(q^{2t-1}+q^{2t+1})\nonumber\\
    =&\frac{(q^{t+1};q^{t+1})_\infty}{(q;q)_\infty}\frac{\sum_{n=1}^{t}{q^{2n}}-q^{t+1}}{1-q^{2t+2}}-2\frac{(q^{t};q^{t})_\infty}{(q;q)_\infty}\frac{\sum_{n=1}^{t-1}{q^{2n}}}{1-q^{2t}}\nonumber\\
    &+\frac{(q^{t};q^{t})_\infty}{(q;q)_\infty}\frac{q^t}{1-q^{2t}}+\frac{(q^{t+1};q^{t+1})_\infty}{(q;q)_\infty}\frac{\sum_{n=1}^{t}{q^{2n}}}{1-q^{2t+2}}\nonumber\\
    &+\frac{(1-q^{t+1})(q^{3t+3};q^{t+1})_\infty}{(q;q)_\infty}(q^{2t+1}+q^{2t+3})-\frac{(1-q^{t})(q^{3t};q^{t})_\infty}{(q;q)_\infty}(q^{2t-1}+q^{2t+1})\nonumber\\
    =&a_{t}-2b_t+c_t+d_{t}+e_{t}-f_{t}.
\end{align}
\qed

We proceed to give the combinatorial interpretations of the six parts respectively. Let $A_{t}(n)$ denote the set of $t+1$-regular OPO-overpartitions $\lambda$ of $n$ with $2\mid o(\lambda)$. 
We further partition the set $A_{t}(n)$ into two subsets $A^1_{t}(n)$ and $A^2_{t}(n)$ as follows:
\begin{itemize}
    \item $A^1_{t}(n)$ is the subset of OPO-overpartitions $\lambda$ in $A_{t}(n)$ such that $f_{\lambda}(1)\ge \sum_{\lambda_i} g(\lambda_i)$;
    \item $A^2_{t}(n)$ is the subset of OPO-overpartitions $\lambda$ in $A_{t}(n)$ such that $f_{\lambda}(1)< \sum_{\lambda_i} g(\lambda_i)$.
\end{itemize}
Here and throughout this paper, for fixed $t$ and let $m=t^k\cdot c$, where $k,c$ are non-negative integers and  $t\nmid c$. We use $g(m)$ to denote the following function:
\begin{equation}
    g(m)=(t+1)^k\cdot c-m.
\end{equation}

We proceed to show that the generating function of $A_{t}(n)$ coincides with $a_{t}.$

\begin{lem}\label{at}
    We have
\begin{equation}
    \sum_{n=0}^\infty A_{t}(n)q^n=a_{t}.
\end{equation}
\end{lem}

\pf Clearly, 
\[\frac{(q^{t+1};q^{t+1})_\infty}{(q;q)_\infty}\]
equals the generating function of $t+1$-regular partitions, which is the generating function of the non-overlined part in $A_{t}(n)$. Moreover,
\[\frac{\sum_{n=1}^{t}{q^{2n}}-q^{t+1}}{1-q^{2t+2}}\]
represents an even number which is not the multiple of $t+1$, this represents the unique overlined part of $A_{t}(n)$. This completes the proof.\qed

Let $B_t(n)$ denote the set of $t$-regular OPO-overpartitions $\lambda$ of $n$ with $2\mid o(\lambda)$. We proceed to show that the generating function of $B_t(n)$ coincides with $b_t$.

\begin{lem}\label{bt}
    We have 
    \begin{equation}
    \sum_{n=0}^{\infty}B_t(n)q^n=b_t.
    \end{equation}
\end{lem}

\pf Clearly,
$$\frac{(q^{t};q^{t})_\infty}{(q;q)_\infty}$$
equals the generating function of $t$-regular partitions, which is the generating function of the non-overlined part in $B_{t}(n)$. Moreover, 
$$\frac{\sum_{n=1}^{t-1}{q^{2n}}}{1-q^{2t}}$$
represents an even number which is not a multiple of $2t$, this represents the unique overlined part of $B_t(n)$. This completes the proof.\qed

Let $C_t(n)$ denote the set of OPO-overpartitions $\lambda$ of $n$. Moreover, we require that $o(\lambda)=(2k-1)t$ for some $k\ge 1$. Furthermore, all the other non-overline parts cannot be multiples of $t$. We proceed to show that the generating function of $C_t(n)$ coincides with $c_t$.
\begin{lem}\label{ct}
    We have 
    \begin{equation}
    \sum_{n=0}^{\infty}C_t(n)q^n=c_t.
    \end{equation}
\end{lem}

\pf Similarly, 
$$\frac{(q^{t};q^{t})_\infty}{(q;q)_\infty}$$
equals the generating function of $t$-regular partitions, which is the generating function of the non-overlined part in $C_{t}(n)$. Moreover, 
$$\frac{q^t}{1-q^{2t}}$$
represents an odd multiples of $t$, and this represents the  overlined part of $C_t(n)$. This completes the proof.
\qed

Let $D_{t}(n)$ denote the set of OPO-overpartitions $\lambda$ of $n$ satisfies the following restrictions:
\begin{itemize}
    \item[(1)] All the non-overlined part cannot be divided by $t+1$;
    \item[(2)] $2|o(\lambda)$ and $2(t+1)\nmid o(\lambda)$.
\end{itemize}
We proceed to show that the generating function of $D_{t}(n)$ coincides with $d_{t}$.

\begin{lem}\label{dt}
    We have 
    \begin{equation}
    \sum_{n=0}^\infty D_{t}(n)q^n=d_{t}.
    \end{equation}
\end{lem}

\pf Similarly, 
$$\frac{(q^{t};q^{t})_\infty}{(q;q)_\infty}$$
equals the generating function of $t$-regular partitions, which is the generating function of the non-overlined part in $D_{t}(n)$. Moreover, 
$$\frac{\sum_{n=1}^{t}{q^{2n}}}{1-q^{2t+2}}$$
represents that the overlined part $o(\lambda)$ satisfies that  $2\mid o(\lambda)$ and $2t+2\nmid o(\lambda)$. This completes the proof.\qed 

Let $E_{t}(n)$ denote the set of OPO-overpartitions $\lambda$ of $n$ satisfies that  $o(\lambda)=2t+1$ or $2t+3$. Moreover, each non-overlined part is either $2t+2$ or cannot be divided by $t+1$. We proceed to show that the generating function of $E_{t}(n)$ coincides with $e_{t}$. 
\begin{lem}\label{et}
    We have
    \begin{equation}
    \sum_{n=0}^\infty E_{t}(n)q^n=e_{t}.
    \end{equation}
\end{lem}

\pf Clearly, 
$$\frac{(1-q^{t+1})(q^{3t+3};q^{t+1})_\infty}{(q;q)_\infty}$$
equals the generating function of partitions, where each part is either $2t+2$ or not divisible by $t+1$. Moreover,
$$(q^{2t+1}+q^{2t+3})$$
represents that the only overlined part is either $2t+1$ or $2t+3$. This completes the proof. \qed

Let $F_t(n)$ denote the set of OPO-overpartitions $\lambda$ of $n$ with   $o(\lambda)=2t-1$ or $2t+1$. Furthermore, each non-overlined part is either $2t$ or cannot be divided by $t$. We proceed to show the generating function $F_t(n)$ coincides with $f_t$.
\begin{lem}\label{ft}
    We have
    \begin{equation}
        \sum_{n=0}^\infty F_t(n)q^n=f_{t}.
    \end{equation}
\end{lem}
\pf Similarly,
$$\frac{(1-q^{t})(q^{3t};q^{t})_\infty}{(q;q)_\infty}$$
equals the generating function of partitions, which need that each part is either $2t$ or cannot be divided by $t$. Moreover,
$$(q^{2t-1}+q^{2t+1})$$
represents that the only overlined part is either $2t-1$ or $2t+1$. This completes the proof. \qed

In order to prove the non-negativity of \eqref{bt+1-bt}, in the rest of this section, we will give three injections $\phi_1,\phi_2$ and $\phi_3$ as shown in the following table:
\[\begin{array}{lllll}
    \phi_1\colon & B_t(n)&\rightarrow & A_t^1(n)\cup C_t(n)\\
    \phi_2\colon & B_t(n)&\rightarrow & D_t(n)\\
    \phi_3\colon & F_t(n)&\rightarrow & A_t^2(n)\cup E_t(n)\\
     &&& 
\end{array}\]
Thus \eqref{bt+1-bt} has non-negative coefficients.

Now, we give the first injection $\phi_1$ from $B_t(n)$ into $A_{t}(n)$ and $C_t(n)$.

\begin{lem}\label{lem-phi_1}
    There is an injection $\phi_1$ from $B_t(n)$ into $A_{t}^1(n)\cup C_t(n)$.
\end{lem}
    
\pf Given $\lambda \in B_t(n)$, by the definition of $B_t(n)$, $\lambda$ is a $t$-regular OPO-overpartitions with  $2\mid o(\lambda)$. Assume $\lambda=(\lambda_1,\ldots,\lambda_\ell)$, where $\overline{o(\lambda)}$ is overlined. We define a map $\alpha$ from $\mathbb{N}$ to the set of overpartitions, satisfies that for any $m\in\mathbb{N}$,
\[m=|\alpha(m)|.\]
There are two cases, based on the overlined part $\overline{o(\lambda)}$.
\begin{itemize}
    \item [Case $1$] $o(\lambda)=c(t+1)^k$,  where $k\ge 0$, $t+1\nmid c$ and $c$ is even, then  set $\alpha(\overline{o(\lambda)})=(\overline{ct^k},1^{g(ct^k)})$. For non-overlined part $\lambda_j$, write $\lambda_j=c'(t+1)^{k'}$, where $k',c'$ are nonnegative integers. Then define $\alpha(\lambda_j)=(c't^{k'},1^{g(c't^{k'})})$.
    \item [Case $2$]  $o(\lambda)=c(t+1)^k$,  where $k\ge 1$, $t+1\nmid c$ and $c$ is odd, set $\alpha(\overline{o(\lambda)})=(\overline{ct^k},1^{g(ct^k)})$. For non-overlined $\lambda_j$, set $\alpha(\lambda_j)=\lambda_j$.
\end{itemize}

Now, we give the injection $\phi_1$ from $B_t(n)$ into $A_{t}^1(n)\cup C_t(n)$ as follows. Define

\begin{equation}\label{eq-def-mu-phi1}
    \mu=\phi_1(\lambda)=\bigcup_{\lambda_i\in\lambda}\alpha(\lambda_i),
\end{equation}
where $\lambda_i$ ranges over all parts of $\lambda$, whether it is overlined or not.

In Case $1$, it is easy to check that the overlined part of $\alpha(\overline{o(\lambda)})$ is even. Moreover, we see that $f_\mu(1)\geq \sum_{\mu_i\in\mu}g(\mu_i)$, since for each $\mu_i>1$, by \eqref{eq-def-mu-phi1}, we see that there exists a unique $k$ such that $\mu_i$ is in $\alpha(\lambda_k)$. Thus, by the construction of $\alpha$, we see that $f_{\alpha(\lambda_k)}(1)=g(\mu_i)$. Moreover, it is clear that when $\mu_i=1$, we have $g(1)=0$. Thus, by \eqref{eq-def-mu-phi1} we deduce that $f_\mu(1)\geq \sum_{\mu_i\in\mu}g(\mu_i)$. Hence $\mu \in A_{t}^1(n)$. Moreover, in Case $2$, we can find that the overlined part is odd multiples of $t$, so it is routine to check that $\mu \in C_t(n)$.

We next show that $\phi_1$ is an injection. To this end, let $I_{\phi_1}(n)$ denote the image set of $\phi$, which is a subset of $A_{t}^1(n) \cup C_t(n)$. We proceed to establish a map $\psi$ from $I_{\phi_1}(n)$ into $B_t(n)$, such that for any $\lambda\in B_t(n)$, 
\begin{equation}\label{phi_1}
    \psi_1(\phi_1(\lambda))=\lambda.
\end{equation}
This implies that $\phi_1$ is an injection.

We now describe $\psi_1$. For any $\mu\in I_{\phi_1}(n)$ with overlined part $\overline{o(\mu)}$, we see that either $\mu\in A^1_{t}(n)$ or $\mu\in C_t(n)$. Moreover, by the definition of $\phi_1$,  it is clear that $f_\mu(1)\geq \sum_{\mu_i\in\mu}g(\mu_i)$. 

Define $\beta(n)=n+g(n)$. 
Now, we define:
\begin{equation}
    \lambda=\psi_1(\mu)=\left(\bigcup_{\mu_i\in\mu}\beta(\mu_i)\right)\setminus (1^{\sum_{\mu_i\in\mu}g(\mu_i)}),
\end{equation}
where set $o(\lambda)=\beta(\overline{o(\mu)})$.
And it is routine to check $\psi$ is an injection and \eqref{phi_1} holds. Thus $\phi_1$ is an injection.\qed

For example, set $t=3$, and $\lambda=(10,\overline{8},7,4,1)$. Applying the injection $\phi_1$, it is easy to check that $o(\lambda)=4\cdot 2$ which is in Case 1. Thus $\mu=\phi_1(\lambda)=(10,7,\overline{6},3,1^4)$ and $\mu\in A_{t}^1(n)$. Moreover, it is routine to check that $\psi_1(\phi_1(\lambda))=\lambda$.

For another example, set $t=3$, and $\lambda=(10,8,7,\overline{4},1)$.  Applying the injection $\phi_1$, it is easy to check that $o(\lambda)=4\cdot 1$ which is in Case 2. Thus $\mu=\phi_1(\lambda)=(10,8,7,\overline{3},1^2)$ and $\mu\in C_{t}(n)$. Moreover, it is routine to check that $\psi_1(\phi_1(\lambda))=\lambda$.

Now, we give the second injection $\phi_2$ from $B_t(n)$ into $D_{t}(n)$.

\begin{lem}\label{lem-phi-2}
    There is an injection $\phi_2$ from $B_t(n)$ into $D_{t}(n)$.
\end{lem}

\pf Given $\lambda\in B_t(n)$, by the definition of $B_t(n)$, $\lambda$ is a $t$-regular OPO-overpartitions with $2\mid o(\lambda)$. Assume $\lambda=(\lambda_1,..,\lambda_\ell)$, where $\overline{o(\lambda)}$ is overlined. We define a map $\gamma$ from $\mathbb{N}$ to the set of overpartitions satisfies that for any $m\in \mathbb{N}$,
$$m=|\gamma(m)|.$$
For the non-overlined part $\lambda_j$, write $\lambda_j=c'(t+1)^{k'}$, where $k',c'$ are nonnegative integers. Then define $\gamma(\lambda_j)=(c't^{k'},1^{g(c't^{k'})})$. For the overlined part $\overline{o(\lambda)}$, assume $o(\lambda)=c(t+1)^k$, where $k\ge 0$ and $t+1\nmid c$. Define
\begin{equation}
    \gamma(\overline{o(\lambda)})=\begin{cases}
        \overline{o(\lambda)}&\text{if $k=1$ and $c$ is odd.}\\
        (\overline{ct^{k-1}(t+1)},1^{(t+1)^{k}-t^{k-1}(t+1)})&\text{if $k\geq 2 $ and $c$ is odd.}\\
        (\overline{ct^k},1^{g(ct^k)})&\text{otherwise.}\\
        
    \end{cases}
\end{equation}

Now, we give the injection $\phi_2$ from $B_t(n)$ into $D_{t}(n)$.
\begin{equation}
    \mu=\phi_2(\lambda)=\bigcup_{\lambda_i\in\lambda}\gamma(\lambda_i).
\end{equation}
where $\lambda_i$ ranges over all the parts of $\lambda$, whether is overlind or not.

We proceed to show that $\mu\in D_t(n)$. When $c$ is even, we see that $o(\mu)=ct^k$ which implies $2t+2\nmid o(\mu)$. Thus $\mu\in D_t(n)$. We next assume that $c$ is odd. From $2\mid o(\lambda)=c(t+1)^k$ we have $k\ge 1$. When $k=1$ we have $o(\mu)=c(t+1)$, thus clearly $2\mid o(\mu)$ and $2(t+1)\nmid o(\mu)$, again we have $\mu\in D_t(n)$. Finally, when $c$ is odd and $k\ge 2$, we have $o(\mu)={ct^{k-1}(t+1)}$, which is an odd multiple of $t+1$. From the above analysis, we conclude that $\mu\in D_t(n)$. 

We next show that $\phi_2$ is an injection. To this end, let $I_{\phi_2}(n)$ denote that the image set of $\phi_2$, which is a subset of $D_{t}(n)$. We proceed to establish a map $\psi_2$ from $I_{\phi_2}(n)$ into $B_t(n)$, such that for any $\lambda\in B_t(n)$, 
\begin{equation}\label{phi_2}
    \psi_2(\phi_2(\lambda))=\lambda.
\end{equation}
This implies $\phi_2$ is an injection.

We now describe $\psi_2$. For any $\mu\in I_{\phi_2}(n)$, it is clear that $2\mid o(\mu)$ and $2(t+1)\nmid o(\mu)$. Moreover, by the definition of $\phi_2$, we see that $f_\mu(1)\geq \sum_{\mu_j\in\mu}g(\mu_j)$.

Define $\delta(n)=n+g(n)$. Now, we define:
\begin{equation}
    \lambda=\psi_2(\mu)=\left(\bigcup_{\mu_j\in\mu}\delta(\mu_j)\right)\setminus(1^{\sum_{\mu_j\in\mu}g(\mu_j)}),
\end{equation}
where $\overline{\delta(o(\mu))}$ is the overlined part. And it is routine to check $\psi_2$ is an injection and \eqref{phi_2} holds. Thus $\phi_2$ is an injection.\qed

For example, set $t=3$, $\lambda=(10,8,7,\overline{4},1)$. Applying the injection of $\phi_2$, it is easy to check that $\mu=\phi_2(\lambda)=(10,7,6,\overline{4},1^3)$. For another example, define $\lambda=(\overline{16},8,4,1,1)$, it is easy to check that $\mu=\phi_2(\lambda)=(\overline{12},6,3,1^9)$. For the third instance, let $\lambda=(10,\overline{8},7,4,1)$, it is easy to check that $\mu=\phi_2(\lambda)=(10,7,\overline{6},3,1^4)$. Moreover, it is routine to check that $\psi_2(\mu)=\lambda$ in all the above three examples.

Now we give the third injection $\phi_3$ from $F_t(n)$ into $E_{t}(n)\cup A_{t+1}^2(n)$.

\begin{lem}\label{lem-phi3}
    There is an injection $\phi_3$ from $F_t(n)$ into $A_{t}^2(n)\cup E_{t}(n)$.
\end{lem}

\pf Given $\lambda \in F_t(n)$, by the definition of $F_t(n)$, $\lambda$ is a OPO-overpartition of $n$ with $o(\lambda)=2t-1$ or $2t+1$. Moreover, each non-overlined part is either $2t$ or cannot be divide by $t$. Assume $\lambda=(\lambda_1,...,\lambda_\ell)$. We define a map $\epsilon$ from $\mathbb{N}$ to the set of overpartitions as follows.

 For the non-overlined part $\lambda_j=c'(t+1)^{k'}$, set $\epsilon(\lambda_j)=(c't^{k'},1^{g(c't^{k'})})$ except $(c',k')=(2,1)$; for the part $2t+2$, set $\epsilon(2t+2)=2t+2$.
We next define $\epsilon(\overline{o(\lambda)})$ and $\mu=\phi_3(\lambda)$. There are $8$ cases, based on the overlined part $\overline{o(\lambda)}$ and the structure of $\lambda$.

\begin{itemize}
    \item [Case $1$] If $o(\lambda)=2t+1$, define $\epsilon(\overline{2t+1})=\overline{2t+1}$ and
    \begin{equation}
        \mu=\phi_3(\lambda)=\bigcup_{\lambda_j\in\lambda}\epsilon(\lambda_j).
    \end{equation}
    \item [Case $2$] If $o(\lambda)=2t-1$, and $f_\lambda(1)\geq 4$. Set $\epsilon(\overline{2t-1})=\overline{2t+3}$ and
    \begin{equation}
        \mu=\phi_3(\lambda)=\bigcup_{\lambda_j\in\lambda}\epsilon(\lambda_j)\setminus(1^4).
    \end{equation}
    \item [Case $3$] If $o(\lambda)=2t-1$, $f_\lambda(1)\leq 3$ and $f_\lambda(2t+2)\geq 2$. Define $\epsilon(\overline{2t-1})=\overline{2t+3}$ and
    \begin{equation}
        \mu=\phi_3(\lambda)=\bigcup_{\lambda_j\in\lambda}\epsilon(\lambda_j)\cup(t,t,t,t)\setminus(2t+2,2t+2).
    \end{equation}
    \item [Case $4$] If $o(\lambda)=2t-1$, $f_\lambda(1)=0$ and $f_\lambda(2t+2)=0$. Set $\epsilon(\overline{2t-1})=(\overline{t-1},t)$ and
    \begin{equation}
        \mu=\phi_3(\lambda)=\bigcup_{\lambda_j\in\lambda}\epsilon(\lambda_j).
    \end{equation}
    \item[Case $5$] If $o(\lambda)=2t-1$, $f_\lambda(1)=1$ or $2$ and $f_\lambda(2t+2)=0$. Let $\epsilon(\overline{2t-1})=\overline{2t}$ and
    \begin{equation}
       \mu=\phi_3(\lambda)=\bigcup_{\lambda_j\in\lambda}\epsilon(\lambda_j)\setminus(1). 
    \end{equation}
    \item [Case $6$] If $o(\lambda)=2t-1$, $f_\lambda(1)=3$ and $f_\lambda(2t+2)=0$. Set $\epsilon(\overline{2t-1})=(t,t,\overline{2})$ and
    \begin{equation}
        \mu=\phi_3(\lambda)=\bigcup_{\lambda_j\in\lambda}\epsilon(\lambda_j)\setminus(1^3).
    \end{equation}

    \item[Case $7$] If $o(\lambda)=2t-1$, $f_\lambda(1)=0$ and $f_\lambda(2t+2)=1$. Set $\epsilon(\overline{2t-1})=(\overline{2t},t,t,1)$ and 
    \begin{equation}
        \mu=\phi_3(\lambda)=\bigcup_{\lambda_j\in\lambda}\epsilon(\lambda_j)\setminus(2t+2).
    \end{equation}
    \item[Case $8$] If $o(\lambda)=2t-1$, $f_\lambda(1)=1$, $2$ or $3$ and $f_\lambda(2t+2)=1$. Set $\epsilon(\overline{2t-1})=(2t,2t,\overline{2})$ and 
    \begin{equation}
    \mu=\phi_3(\lambda)=\bigcup_{\lambda_j\in\lambda}\epsilon(\lambda_j)\setminus(2t+2,1).
    \end{equation}
\end{itemize} 

We first show that the image sets of the first three cases are in $E_t(n)$. It is clear that in the first three cases $o(\mu)=2t+1$ or $2t+3$, and each non-overlined part is either $2t+2$ or cannot be divided by $t+1$. This implies the image sets of the first three cases are in $E_t(n)$. We proceed to show that the image sets of the last five cases are in $A_t^2(n)$. It is trivial to check that $o(\mu)=2,t-1,2t$, which is even number and is not multiple of $t+1$. Moreover, from the construction of $\phi_3$, it is easy to check that $f_\mu(2t+2)=0$, which implies $\mu$ is $t+1$-regular partition in all the five cases. Furthermore, we calculate $f_{\mu}(1)- \sum_{\mu_i} g(\mu_i)$ in each case as listed in Table \ref{tb-1}. Thus we show that the image sets of the last five cases are in $A_t^2(n)$. Moreover, from Table \ref{tab:my_label}, we see that the image sets of the eight cases are pairwise non-intersect. 
\begin{table}[]
    \centering
    \begin{tabular}{c|c}
    \hline
    & $f_{\mu}(1)- \sum_{\mu_i} g(\mu_i)$\\
    \hline
    Case 4 & $\leq-1$\\
    \hline
    Case 5 & $\leq-1$\\
    \hline
    Case 6 & $\leq-2$\\
    \hline
    Case 7 & $\leq-3$\\
    \hline
    Case 8 & $\leq-2$\\
    \hline
    \end{tabular}
    \caption{List of $f_{\mu}(1)- \sum_{\mu_i} g(\mu_i)$ in Case 4 $\sim$ Case 8.}
    \label{tb-1}
\end{table}

\begin{table}[]
    \centering
    \begin{tabular}{c|c}
    \hline
    & $E_{t}(n)$\\
    \hline
    Case 1     &  $o(\mu)=2t+1$\\
    \hline
    Case 2 &  $o(\mu)=2t-1$, $ f_\mu(1)\geq \sum_{\mu_j\not=2t}g(\mu_j)$\\
    \hline
    Case 3 &       $o(\mu)=2t-1$, $ f_\mu(1)< \sum_{\mu_j\not=2t}g(\mu_j)$  \\
    \hline
    & $A_{t}^2(n)$\\
    \hline
    Case 4 & $o(\mu)=t-1$, $f_\mu(1)-\sum_{\mu_j\not=2t}g(\mu_j)=-1$\\
    \hline
    Case 5 & $o(\mu)=2t$, $f_\mu(1)-\sum_{\mu_j\not=2t}g(\mu_j)=0$ or $1$\\
    \hline
    Case 6 & $o(\mu)=2$, $f_\mu(1)-\sum_{\mu_j\not=2t}g(\mu_j)=-2$\\
    \hline
    Case 7 & $o(\mu)=2t$, $f_\mu(1)-\sum_{\mu_j\not=2t}g(\mu_j)=-1$\\
    \hline
    Case 8 & $o(\mu)=2$, $f_\mu(1)-\sum_{\mu_j\not=2t}g(\mu_j)=0$, $1$ or $2$\\
    \hline
    
    \end{tabular}
    \caption{The image sets of eight cases in Lemma \ref{lem-phi3}}
    \label{tab:my_label}
    
\end{table}

We proceed to show that $\phi_3$ is an injection. To this end, let $I_{\phi_3}(n)$ denote the image set of $\phi_3$, which is a subset of $A_{t}^2(n) \cup E_{t}(n)$. We proceed to establish a map $\psi_3$ from $I_{\phi_3}(n)$ into $F_t(n)$, such that for any $\lambda\in F_t(n)$, 
\begin{equation}\label{phi_3}
    \psi_3(\phi_3(\lambda))=\lambda.
\end{equation}
This implies $\phi_3$ is an injection.

We now describe $\psi_3$. For any partition $\mu\in I_{\phi_3}(n)$, we first define a map $\zeta$. For the non-overlined part $\mu_j$, define
\begin{equation}\label{equ-def-zeta}
    \zeta(\mu_j)=\begin{cases}
        \mu_j+g(\mu_j)&\text{if } \mu_j\ne 2t;\\
        \mu_j&\text{if }\mu_j=2t.
    \end{cases}
\end{equation}
For the overlined part $\overline{o(\mu)}$, if  $\mu\in  E_{t}(n)$, there are three cases.
\begin{itemize}
    \item [Case $1$] If $o(\mu)=2t+1$, by definition, it is easy to check that $ f_\mu(1)\geq \sum_{\mu_j\ne 2t}g(\mu_j)$. Define $\zeta(\overline{2t+1})=\overline{2t+1}$ and
    \begin{equation}
        \lambda=\psi_3(\mu)=\bigcup_{\mu_j\in\mu}\zeta(\mu_j)\setminus(1^{\sum_{\mu_j\ne 2t}g(\mu_j)}).
    \end{equation}
    \item [Case $2$] If $o(\mu)=2t+3$ and $ f_\mu(1)\geq \sum_{\mu_j\ne 2t}g(\mu_j)$. Define $\zeta(\overline{2t+3})=(\overline{2t-1},1^4)$ and 
    \begin{equation}
    \lambda=\psi_3(\mu)=\bigcup_{\mu_j\in\mu}\zeta(\mu_j)\setminus(1^{\sum_{\mu_j\ne 2t}g(\mu_j)}).
    \end{equation}
    \item [Case $3$] If $o(\mu)=2t+3$ and $ f_\mu(1)< \sum_{\mu_j\ne 2t}g(\mu_j)$. By Table \ref{tab:my_label}, we see that this case is in Case 3, which implies $f_\mu(t)\ge 4$. Set $\nu=\mu\setminus(t^4)$. Define $\zeta(\overline{2t+3})=((2t+2)^2,\overline{2t-1})$ and 
    \begin{equation}
         \lambda=\psi_3(\mu)=\bigcup_{\nu_j\in\nu}\zeta(\nu_j)\setminus(1^{\sum_{\nu_j\ne 2t}g(\nu_j)}).
    \end{equation}
\end{itemize}

If $\mu\in  A_{t}^2(n)$,  then there are five cases for  $\zeta(\overline{o(\mu)})$. 

\begin{itemize}
    \item [Case $1$] If $o(\mu)=t-1$. Then by Table \ref{tab:my_label}, we have $\mu$ is in the image set of Case $4$; and by the construction of $\phi_3$, we see that $f_\mu(t)\geq1$. Set $\nu=\mu\setminus(t)$, and $\zeta(\overline{t-1})=\overline{2t-1}$. Define,
    \begin{equation}
        \lambda=\psi_3(\mu)=\bigcup_{\nu_j\in\nu}\zeta(\nu_j)\setminus(1^{\sum_{\nu_j\ne 2t}g(\nu_j)}).
    \end{equation}
    \item [Case $2$] If $o(\mu)=2t$, and $f_\mu(1)-\sum_{\mu_j\not=2t}g(\mu_j)=0$ or $1$. In this case, by Table \ref{tab:my_label}, we have $\mu$ is in the image set of Case $5$. Set $\zeta(\overline{2t})=(\overline{2t-1},1)$. Define, 
    \begin{equation}
    \lambda=\psi_3(\mu)=\bigcup_{\mu_j\in\mu}\zeta(\mu_j)\setminus(1^{\sum_{\mu_j\ne 2t}g(\mu_j)}).
    \end{equation}
    \item [Case $3$] If $o(\mu)=2$ and $f_\mu(1)-\sum_{\mu_j\not=2t}g(\mu_j)=-2$. In this case,  by Table \ref{tab:my_label}, we have $\mu$ is in the image set of Case $6$. Thus by the construction of $\phi_3$, we see that $f_{\mu}(t)\ge 2$. Set $\nu=\mu\setminus(t,t)$ and $\zeta(\overline{2})=(\overline{2t-1},1^3)$. Define, 
    \begin{equation}
        \lambda=\psi_3(\mu)=\bigcup_{\nu_j\in\nu}\zeta(\nu_j)\setminus(1^{\sum_{\nu_j\ne 2t}g(\nu_j)}).
    \end{equation}
    \item [Case $4$] If $o(\mu)=2t$, $f_\mu(1)-\sum_{\mu_j\not=2t}g(\mu_j)=-1$. In this case,  by Table \ref{tab:my_label}, we have $\mu$ is in the image set of Case $7$. Thus by the construction of $\phi_3$, we see that $f_\mu(t)\geq 2$ and $f_\mu(1)\geq 1$. Set $\nu=\mu\setminus(t,t,1)$ and $\zeta(\overline{2t})=(\overline{2t-1},2t+2)$. Define,
    \begin{equation}
        \lambda=\psi_3(\mu)=\bigcup_{\nu_j\in\nu}\zeta(\nu_j)\setminus(1^{\sum_{\nu_j\ne 2t}g(\nu_j)}).
    \end{equation}
    \item[Case $5$] The overlined part $o(\mu)=2$, $f_\mu(1)-\sum_{\mu_j\not=2t}g(\mu_j)=0$, $1$ or $2$. In this case,  by Table \ref{tab:my_label}, we have $\mu$ is in the image set of Case $8$. Thus by the construction of $\phi_3$, we see that $f_\mu(2t)\geq 2$. Set $\nu=\mu\setminus(2t,2t)$ and $\zeta(\overline{2})=(\overline{2t-1},2t+2,1)$. Define,
    \begin{equation}
         \lambda=\psi_3(\mu)=\bigcup_{\nu_j\in\nu}\zeta(\nu_j)\setminus(1^{\sum_{\nu_j\ne 2t}g(\nu_j)}).
    \end{equation}
\end{itemize}

It is routine to check $\psi_3$ is an injection and \eqref{phi_3} holds. Thus $\phi_3$ is an injection. \qed

For example, set $t=3$, $\lambda^1=(8,\overline{7},6,4,2,1)$. Applying the injection $\phi_3$, it is easy to check that $\mu^1=\phi_3(\lambda^1)=(8,\overline{7},6,3,2,1,1)$. For another example, let $t=3$ and $\lambda^2=(8,6,\overline{5},4,1^5)$, it is easy to check that $\mu^2=\phi_3(\lambda^2)=(\overline{9},8,6,3,1^2)$; $\lambda^3=(7,6,\overline{5},2)$, it is easy to check that $\mu^3=\phi_3(\lambda^3)=(7,6,3,\overline{2},2)$. Moreover, it is routine to check $\psi_3(\mu)=\lambda$.

{\noindent \it Proof of Theorem \ref{thm-com-int-umn} for odd number $t$.} From the Lemma \ref{lem-phi_1} $\sim$ \ref{lem-phi3}, it is clear that \eqref{bt+1-bt} has non-negative coefficient of $q^n$. Thus  it proves Theorem \ref{thm-com-int-umn} for odd number $t$.\qed

\section{A proof of Theorem \ref{thm-com-int-umn} for even number $t$}\label{3}

This section is devoted to proving Theorem \ref{thm-com-int-umn} for even number $t$.  Similar as in Section \ref{2}, we first divide the generating function of $b_{t+1,2}(n)-b_{t,2}(n)$ into six parts. We then give the combinatorial interpretations of these six parts respectively. Three injections will be built and Theorem \ref{thm-com-int-umn} follows from these three injections. It should be noted that  case $t=4$ needs to slightly modify the third injection. In this way we show that Theorem \ref{thm-com-int-umn} holds for all $t\ge 4$ even.

Similarly, we first transform the generating function of $b_{t+1,2}(n)-b_{t,2}(n)$. To this end, we define $a_t$, $b_t$, $c_t$, $d_t$, $e_t$ and $f_t$ as follows.
\begin{align*}
    a_t&=\frac{(q^{t+1};q^{t+1})_\infty}{(q;q)_\infty}\frac{\sum_{n=1}^{t}{q^{2n}}}{1-q^{2t+2}},\\
    b_t&=\frac{(q^{t+1};q^{t+1})_\infty}{(q;q)_\infty}\frac{q^{t+1}}{1-q^{2t+2}},\\
    c_t&=\frac{(q^{t};q^{t})_\infty}{(q;q)_\infty}\frac{\sum_{n=1}^{t-1}{q^{2n}-q^t}}{1-q^{2t}},\\
    d_t&=\frac{(q^{t};q^{t})_\infty}{(q;q)_\infty}\frac{\sum_{n=1}^{t-1}{q^{2n}}}{1-q^{2t}},\\
    e_t&=\frac{(1-q^{t+1})(q^{3t+3};q^{t+1})_\infty}{(q;q)_\infty}(q^{2t+1}+q^{2t+3}),\\
    f_t&=\frac{(1-q^{t})(q^{3t};q^{t})_\infty}{(q;q)_\infty}(q^{2t-1}+q^{2t+1}).
\end{align*}
It should be noted that these definitions are different as in Section 2. This does not lead to a ambiguous since $t$ is odd in Section 2 and now we consider the case $t$ is even.

\begin{lem}\label{bt+1-bt1}
    For given even integer $t$, we have that 
    \begin{equation}\label{bt+12-bt2}
         \sum_{n=1}^{\infty}b_{t+1,2}(n)q^n-\sum_{n=1}^{\infty}b_{t,2}(n)q^n=2a_t-b_t-c_t-d_t+e_t-f_t.
    \end{equation}
\end{lem}

\pf By Theorem \ref{bt2}, we see that 
\begin{align}
    &\sum_{n=1}^\infty b_{t+1,2}(n)-\sum_{n=1}^\infty b_{t,2}(n)\nonumber\\
    =&\frac{(q^{t+1};q^{t+1})_\infty}{(q;q)_\infty}\left(\frac{2q^2}{1-q^2}-\frac{q^{t+1}}{1-q^{t+1}}+\frac{q^{2t+1}-q^{2t+2}+q^{2t+3}}{1-q^{2t+2}}\right)\nonumber\\
    &-\frac{(q^t;q^t)_\infty}{(q;q)_\infty}\left(\frac{2q^2}{1-q^2}-\frac{q^t}{1-q^t}+\frac{q^{2t-1}-q^{2t}+q^{2t+1}}{1-q^{2t}}\right)\nonumber\\
    =&\frac{(q^{t+1};q^{t+1})_\infty}{(q;q)_\infty}\left(\frac{2q^2}{1-q^2}-\frac{q^{t+1}}{1-q^{t+1}}-\frac{q^{2t+2}}{1-q^{2t+2}}\right)-\frac{(q^t;q^t)_\infty}{(q;q)_\infty}\left(\frac{2q^2}{1-q^2}-\frac{q^t}{1-q^t}-\frac{q^{2t}}{1-q^{2t}}\right)\nonumber\\
    &+\frac{(q^{t+1};q^{t+1})_\infty}{(q;q)_\infty}\left(\frac{q^{2t+1}+q^{2t+3}}{1-q^{2t+2}}\right)-\frac{(q^t;q^t)_\infty}{(q;q)_\infty}\left(\frac{q^{2t-1}+q^{2t+1}}{1-q^{2t}}\right)\nonumber\\
    =&2\frac{(q^{t+1};q^{t+1})_\infty}{(q;q)_\infty}\left(\frac{q^2}{1-q^2}-\frac{q^{2t+2}}{1-q^{2t+2}}\right)-\frac{(q^{t+1};q^{t+1})_\infty}{(q;q)_\infty}\left(\frac{q^{t+1}}{1-q^{t+1}}-\frac{q^{2t+2}}{1-q^{2t+2}}\right)\nonumber\\
    &-\frac{(q^t;q^t)_\infty}{(q;q)_\infty}\left(\frac{q^2}{1-q^2}-\frac{q^t}{1-q^t}\right)-\frac{(q^t;q^t)_\infty}{(q;q)_\infty}\left(\frac{q^2}{1-q^2}-\frac{q^{2t}}{1-q^{2t}}\right)\nonumber\\
    &+\frac{(1-q^{t+1})(q^{3t+3};q^{t+1})_\infty}{(q;q)_\infty}(q^{2t+1}+q^{2t+3})-\frac{(1-q^{t})(q^{3t};q^{t})_\infty}{(q;q)_\infty}(q^{2t-1}+q^{2t+1})\nonumber\\
    =&2\frac{(q^{t+1};q^{t+1})_\infty}{(q;q)_\infty}\frac{\sum_{n=1}^{t}{q^{2n}}}{1-q^{2t+2}}-\frac{(q^{t+1};q^{t+1})_\infty}{(q;q)_\infty}\frac{q^{t+1}}{1-q^{2t+2}}\nonumber\\
    &-\frac{(q^{t};q^{t})_\infty}{(q;q)_\infty}\frac{\sum_{n=1}^{t-1}{q^{2n}-q^t}}{1-q^{2t}}-\frac{(q^{t};q^{t})_\infty}{(q;q)_\infty}\frac{\sum_{n=1}^{t-1}{q^{2n}}}{1-q^{2t}}+\frac{(1-q^{t+1})(q^{3t+3};q^{t+1})_\infty}{(q;q)_\infty}(q^{2t+1}+q^{2t+3})\nonumber\\
    &-\frac{(1-q^{t})(q^{3t};q^{t})_\infty}{(q;q)_\infty}(q^{2t-1}+q^{2t+1})\nonumber\\
    =&2a_t-b_t-c_t-d_t+e_t-f_t.
\end{align}
\qed

We proceed to give the combinatorial interpretations of the six parts respectively. 
\begin{itemize}
    \item Let $A_t(n)$ denote the set of the $t+1$-regular OPO-overpartitions $\lambda$ of $n$ with $o(\lambda)$ even. Since the coefficient of $a_t$ in \eqref{bt+12-bt2} is $2$, we introduce $\hat{A}_t(n)$  as a copy of $A_t(n)$ with the overlined part tagged. For example, for $t=4$, $(7,\overline{6},6,2,1)\in A_4(22)$ and $(7,\overline{6}',6,2,1)\in \hat{A}_4(22)$. We further partition  ${A}_t(n)$  into two subsets ${A}_t^1(n)$ and ${A}_t^2(n)$, while $\hat{A}_t(n)$ is partitioned into two subsets 
 $\hat{A}_t^3(n)$ and $\hat{A}_t^4(n)$. The specific definitions of these four subsets are given as follows:  
\begin{itemize}
    \item[(1)] $A^1_{t}(n)$ is the subset of OPO-overpartitions $\lambda$ in $A_{t}(n)$ such that $f_{\lambda}(1)\ge \sum_{\lambda_i} g(\lambda_i)$ for $o(\lambda)=ct^k$ when $2\mid c$, and $f_{\lambda}(1)\ge \sum_{\lambda_i} g(\lambda_i)-c(t+1)^{k-1}$ for $o(\lambda)=ct^k$ when $2\nmid c$.
    \item[(2)] $A^2_{t}(n)$ is the subset of OPO-overpartitions $\lambda$ in $A_{t}(n)$ such that $f_{\lambda}(1)< \sum_{\lambda_i} g(\lambda_i)$ for $o(\lambda)=ct^k$ when $2\mid c$, and $f_{\lambda}(1)< \sum_{\lambda_i} g(\lambda_i)-c(t+1)^{k-1}$ for $o(\lambda)=ct^k$ when $2\nmid c$.
    \item[(3)] $\hat{A}^3_t(n)$ is the subset of tagged OPO-overpartitions $\lambda$ in $\hat{A}_t(n)$ such that $f_\lambda(1)\geq \sum_{\lambda_i}g(\lambda_i)$ for $o(\lambda)=ct^k$ when $2\mid c$ and $k\geq 0$; and $f_\lambda(1)\geq g(o(\mu))$ for $o(\lambda)=ct^k$ when $2\nmid c$ and $k\geq 1$.
    \item[(4)] $\hat{A}^4_t(n)$ is the subset of tagged OPO-overpartitions $\lambda$ in $\hat{A}_t(n)$ such that $f_\lambda(1)< \sum_{\lambda_i}g(\lambda_i)$ for $o(\lambda)=ct^k$ when $2\mid c$ and $k\geq 0$; and $f_\lambda(1)< g(o(\mu))$ for $o(\lambda)=ct^k$ when $2\nmid c$ and $k\geq 1$.
\end{itemize}
\item Let $B_t(n)$ denote the set of OPO-overpartitions $\lambda$. Moreover, we require that $o(\lambda)$ is the odd multiple of $t+1$. Furthermore, all the other non-overlined parts cannot be multiples of $t+1$. 
\item Let $C_t(n)$ denote the set of the $t$-regular OPO-overpartitions $\lambda$ of $n$ with $o(\lambda)$ even.
\item Let $D_t(n)$ denote the set of OPO-overpartitions $\lambda$ of $n$ satisfies the following restrictions:
\begin{itemize}
    \item [(1)] All the non-overlined part cannot be divide by $t$;
    \item [(2)]   $2\mid  o(\lambda)$ and  $2t\nmid o(\lambda)$.
\end{itemize}
\item Let $E_{t}(n)$ denote the set of OPO-overpartitions $\lambda$ of $n$ with $o(\lambda)$ equals either  $2t+1$ or $2t+3$. Moreover, each non-overlined part is either $2t+2$ or cannot be divide by $t+1$.
\item Let $F_t(n)$ denote the set of OPO-overpartitions $\lambda$ of $n$ with $o(\lambda)$ equals either  $2t-1$ or $2t+1$. Moreover, each non-overlined part is either $2t$ or cannot be divided by $t$.
\end{itemize}
Using the same argument as in Lemma \ref{at} $\sim$ \ref{ft}, we see that $a_t$, $b_t$, $c_t$, $d_t$, $e_t$ and $f_t$ are the generating functions of $A_t(n)$, $B_t(n)$, $C_t(n)$, $D_t(n)$, $E_t(n)$ and $F_t(n)$ respectively. We omit the detailed proofs.

In order to prove the non-negativity of \eqref{bt+1-bt1},  we will give three injections $\zeta_1,\zeta_2$ and $\zeta_3$ as shown in the following table:
\[\begin{array}{lllll}
    \zeta_1\colon & B_t(n)\cup C_t(n)&\rightarrow & \hat{A}^3_t(n)\\
    \zeta_2\colon & D_t(n)&\rightarrow & A_t^1(n)\\
    \zeta_3\colon & F_t(n)&\rightarrow & A_t^2(n)\cup \hat{A}_t^4(n)\cup E_t(n)\\
     &&& 
\end{array}\]
Thus \eqref{bt+1-bt1} has non-negative coefficients.

We describe the first injection $\zeta_1$ from $B_t(n)\cup C_t(n)$ into $\hat{A}_t(n)$ as follows.

\begin{lem}\label{zeta1}
    When $t\ge 4$ even, there is an injection $\zeta_1$ from $B_t(n)\cup C_t(n)$ into $\hat{A}_t^3(n)$.
\end{lem}

\pf On the one hand, given $\lambda \in B_t(n)$, by definition $\lambda$ is an OPO-overpartition with $o(\lambda)$ is odd multiple of $t+1$. Moreover, the non-overlined part is not the multiple of $t+1$.   We define a map $\alpha_1$ from $\mathbb{N}$ to the set of  overpartitions, satisfies that for any $m\in \mathbb{N}$,
$$m=|\alpha_1(m)|.$$

For the overlined part $o(\lambda)=c(t+1)^k$ where $2\nmid c$, $t+1\nmid c$ and $ k\geq1$, define
$$\alpha_1(\overline{o(\lambda)})=(\overline{ct^{k}}',1^{g(ct^{k})}).$$

For the non-overlined part $\lambda_j$, define 
$$\alpha_1(\lambda_j)=\lambda_j.$$
Now we define $\zeta_1(\lambda)$ as follows:
\begin{equation}
    \mu=\zeta_1(\lambda)=\bigcup_{\lambda_i\in\lambda}\alpha_1(\lambda_i).
\end{equation}

On the other hand, for any $\lambda \in C_t(n)$, by the definition of $C_t(n)$, $\lambda$ is a $t$-regular OPO-overpartition with $o(\lambda)$ even.  We define a map $\alpha_1'$ from $\mathbb{N}$ to the set of overpartitions, satisfies that for any $m\in \mathbb{N}$,
$$m=|\alpha_1'(m)|.$$

For the part $o(\lambda)=c(t+1)^k$ where $2\mid c$ and $ k\geq0$, define $\alpha_1'(\overline{o(\lambda)})=(\overline{ct^k}',1^{g(ct^k)}).$; for the non-overlined part $\lambda_j=c'(t+1)^{k'}$, define
$\alpha_1'(\lambda_j)=(c't^{k'},1^{g(c't^{k'})}).$

For the partition $\lambda$ of $C_t(n)$, define 
\begin{equation}
    \mu=\zeta_1(\lambda)=\bigcup_{\lambda_i\in \lambda}\alpha_1'(\lambda_i).
\end{equation}

By the definition of $\alpha_1$ and $\alpha_1'$, for any $\lambda \in B_t(n)$, it is clear that $\mu$ is an OPO-overpartition with $o(\mu)=ct^k$ satisfies that $2\nmid c$, $k\geq1$ and $t\nmid c$, which $f_\lambda(1)\geq g(o(\mu))$; for $\lambda \in C_t(n)$, it can be checked that $\mu$ is an OPO-overpartition with $o(\mu)=ct^k$ such that $2\mid c$ and $t\nmid c$, which $f_\lambda(1)\geq \sum_{\mu_i} g(\mu_i)$. This implies that the image sets of $B_t(n)$ and $C_t(n)$ are non-intersect. Moreover, for the partition $\lambda \in B_t(n)\cup C_t(n)$, it is routine to see that $\mu$ is a $t+1$-regular OPO-overpartition with $o(\mu)$ even. It is clear to check that $\mu \in \hat{A}_t^3(n)$. 

We next show that $\zeta_1$ is an injection. To this end, let $I_{\zeta_1}(n)$ denote that the image set of $\zeta_1$, which is a subset of $\hat{A}_t^3(n)$. We proceed to establish a map $\eta_1$ from $I_{\zeta_1}(n)$ into $B_t(n) \cup C_t(n)$, such that for any $\lambda\in B_t(n) \cup C_t(n)$, 
\begin{equation}\label{eta1}
    \eta_1(\zeta_1(\lambda))=\lambda. 
\end{equation}
This implies $\zeta_1$ is an injection.

We now describe $\eta_1$. Given $\mu\in I_{\zeta_1}(n)$, by definition $\mu$ is a t+1-regular OPO-overpartition with $2\mid o(\mu)$. Set $o(\mu)=ct^k$, where $t\nmid c$,
there are two cases based on either $2\mid c$. 
\begin{itemize}
    \item [Case 1.] If $2\nmid c$, then from the construction of $\zeta_1$, it is clear that $\mu$ is in the image set of $B_t(n)$. Thus we see that $f_\mu(1)\ge g(o(\mu))$. Now set $\beta_1(\overline{o(\mu)}')=\overline{d}$, where $d=o(\mu)+g(o(\mu))$. Moreover, for the non-overlined part $\mu_j$, set $\beta_1(\mu_j)=\mu_j$. Define,
    \begin{equation}
        \eta_1(\mu)=\bigcup_{\mu_i\in \mu}\beta_1(\mu_i)\setminus(1^{g(o(\mu))}).
    \end{equation}
    \item [Case 2.] If $2\mid c$, then from the construction of $\zeta_1$, it is clear that $\mu$ is in the image set of $C_t(n)$. Thus we see that $f_\mu(1)\geq \sum_{\mu_i\in\mu}g(\mu_i)$. Set $\beta_1(\overline{o(\mu)}')=\overline{d}$, where $d=o(\mu)+g(o(\mu))$. For the non-overlined part $\mu_j$, set $\beta_1(\mu_j)=\mu_j+g(\mu_j)$. Define,
    \begin{equation}
        \eta_1(\mu)=\bigcup_{\mu_i\in \mu}\beta_1(\mu_i)\setminus(1^{\sum_{\mu_i\in\mu}g(\mu_i)}).
    \end{equation}
    
\end{itemize}

It is routine to check $\eta_1$ is an injection and \eqref{eta1} holds. Thus $\zeta_1$ is an injection.\qed

For example, set $t=4$ and $\lambda=(11,6,\overline{5},2,1,1)$. Applying the injection $\zeta_1$, it is easy to check that $\lambda$ is in $B_t(n)$. Thus $\mu=\zeta_1(\lambda)=(11,6,\overline{4}',2,1,1,1)$ and $\mu\in \hat{A}_t^3(n)$. Moreover, it is routine to check that $\eta_1(\zeta_1(\lambda))=\lambda$.

For another example, set $t=4$ and $\lambda=(11,\overline{10},6,5,2,1)$. Applying the injection $\zeta_1$, it is easy to check that $\lambda$ is in $C_t(n)$. Thus $\mu=\zeta_2(\lambda)=(11,\overline{8}',6,4,2,1^4)$ and $\mu\in \hat{A}_t^3(n)$. Moreover, it is routine to check that $\eta_1(\zeta_1(\lambda))=\lambda$.

Now, we give the second injection $\zeta_2$ from $D_t(n)$ into $A_t^1(n)$.

\begin{lem}\label{zeta2}
    When $t\ge 4$ even, there is a bijection $\zeta_2$ from $D_t(n)$ into $A_t^1(n)$.
\end{lem}

\pf Given $\lambda \in D_t(n)$, $\lambda$ is the OPO-overpartition with $2\mid o(\lambda)$ and $2t\nmid o(\lambda)$, and the non-overlined part cannot be divided by $t$. Assume $\lambda=(\lambda_1,\lambda_2,...,\lambda_\ell)$, where $\overline{o(\lambda)}$ is overlined. We define a map $\alpha_2$ from $\mathbb{N}$ to the set of overpartitions, satisfies that for any $m \in \mathbb{N}$, 
$$m=|\alpha_2(m)|.$$

For the non-overlined part $\lambda_j=c'(t+1)^{k'}$ where $c'\geq 0$ and $k' \geq 0$, define
$$\alpha_2(\lambda_j)=(c't^{k'},1^{g(c't^{k'})}).$$

We now  construct $\alpha_2(\overline{o(\lambda)})$. Let $o(\lambda)=c(t+1)^k$ where $2\mid c$ and $k \geq 0$, there are two cases.
\begin{itemize}
    \item Case 1. $t\mid c$. By definition we see that $2t\nmid c$, which implies that $\tilde{c}=c/t$ is odd.  Define $$\alpha_2(\overline{o(\lambda)})=(\overline{\tilde{c}t^{k+1}},1^{\tilde{c}t(t+1)^k-\tilde{c}t^{k+1}}).$$
    \item Case 2. $t\nmid c$.  Define $\alpha_2(\overline{o(\lambda)})=(\overline{ct^k},1^{g(ct^k)})$.
\end{itemize}

For the overpartition $\lambda$ of $D_t(n)$, define
\begin{equation}
    \mu=\zeta_2(\lambda)=\bigcup_{\lambda_i\in \lambda}\alpha_2(\lambda_i).
\end{equation}

By the definition of $\alpha_2$, it is clear that $\mu$ is a $t+1$-regular overpartition with the only one even part overlined. Moreover, in Case 1, we have $o(\mu)=ct^k$ with $c$ odd. It is clear that $f_{\lambda}(1)- \sum_{\lambda_i} g(\lambda_i)+c(t+1)^{k-1}\geq 0$. In Case 2, we have $o(\mu)=ct^k$ with $c$ even. Thus  $f_\mu(1)-\sum_{\mu_i\in\mu}g(\mu_i)\geq 0$. This implies  $\mu\in A_t^1(n)$.

We next show that $\zeta_2$ is a bijection. To this end,   we will establish a map $\zeta_2^{-1}$ from $A_t^1(n)$ to $D_t(n)$. For any $\mu \in A_t^1(n)$, by definition we see that $2\mid o(\mu)$ and $2t\nmid o(\mu)$. Moreover, let $o(\mu)=ct^k$, where $t\nmid c$. We have $f_{\mu}(1)- \sum_{\mu_i} g(\mu_i)+c(t+1)^{k-1}\geq 0$ when $2\nmid c$, and $f_{\mu}(1)- \sum_{\mu_i} g(\mu_i)+c(t+1)^{k-1}\geq 0$ when $2\mid c$.
There are two cases.

\begin{itemize}
    \item [Case 1.] $2\mid c$, set $\beta_2(\overline{o(\mu)})=\overline{d},$ 
            where $d=o(\mu)+g(o(\mu))$. And for the non-overlined part $\mu_j$, set 
            $\beta_2(\mu_j)=\mu_j+g(\mu_j)$. Define
            \begin{equation}
            \zeta_2^{-1}(\mu)=\bigcup_{\mu_j\in\mu}(\beta_2(\mu_j))\setminus(1^{\sum_{\mu_i\in\mu}g(\mu_i)}).
            \end{equation}
    \item [Case 2.] $2\nmid c$, set $\beta_2(\overline{o(\mu)})=\overline{d},$ 
            where $d=o(\mu)+ct(t+1)^{k-1}-ct^k=o(\mu)+g(o(\mu))-c(t+1)^{k-1}$. And for the non-overlined part $\mu_j$, set $\beta_2(\mu_j)=\mu_j+g(\mu_j)$. Define
            \begin{equation}
            \zeta_2^{-1}(\mu)=\bigcup_{\mu_j\in\mu}(\beta_2(\mu_j))\setminus(1^{\sum_{\mu_i\in\mu}g(\mu_i)-c(t+1)^{k-1}}).
            \end{equation}
\end{itemize}

It is routine to check $\zeta_2^{-1}$ is the inverse map of $\zeta_2$. \qed

For example, set $t=6$ and $\lambda=(11,7,\overline{6},5,3,2,1)$. Applying the injection $\zeta_2$, it is easy to check that $\lambda$ is in $D_t(n)$. Thus $\mu=\zeta_2(\lambda)=(11,\overline{6},6,5,3,2,1,1)$ and $\mu\in A_t^1(n)$. Moreover, it is routine to check that $\zeta_2^{-1}(\zeta_2(\lambda))=\lambda$.

For another example, set $t=6$ and $\lambda=(14,11,7,6,5,\overline{2},1)$. Applying the injection $\zeta_2$, it is easy to check that $\lambda$ is in $D_t(n)$. Thus $\mu=\zeta_2(\lambda)=(12,11,6,6,5,\overline{2},1^4)$ and $\mu\in A_t^1(n)$. Moreover, it is routine to check that $\zeta_2^{-1}(\zeta_2(\lambda))=\lambda$.

Now, we give the third injection $\zeta_3$ from $F_t(n)$ into $E_t(n)\cup A_t^2(n)\cup \hat{A}_t^4(n)$.

\begin{lem}\label{zeta3}\label{zeta3}
    When $t\ge 4$ even, there is an injection $\zeta_3$ from $F_t(n)$ into $E_t(n)\cup A_t^2(n)\cup \hat{A}_t^4(n)$.
\end{lem}

\pf We first consider the case $t\ge 6$.

Given $\lambda\in F_t(n)$, $\lambda$ is the OPO-overpartition and $o(\lambda)$ is either $2t-1$ or $2t+1$. Furthermore, each non-overlined part of $\lambda$ is either $2t$ or cannot be divided by $t$.   We next define a map $\alpha_3$ from $\mathbb{N}$ to the set of overpartitions.

For the non-overlined part $\lambda_j=c'(t+1)^{k'}$, set $\alpha_3(\lambda_j)=(c't^{k'},1^{g(c't^{k'})})$ except $(c',k')=(2,1)$; for the part $2t+2$, set $\alpha_3(2t+2)=2t+2$.
We next define $\alpha_3(\overline{o(\lambda)})$ and $\mu=\zeta_3(\lambda)$. There are $8$ cases, based on the overlined part $\overline{o(\lambda)}$ and the structure of $\lambda$. 

\begin{itemize}
    \item [Case 1.] $o(\lambda)=2t+1$. In this case, set $\alpha_3(\overline{2t+1})=\overline{2t+1}$. Define
    \begin{equation}
        \mu=\zeta_3(\lambda)=\bigcup_{\lambda_i\in \lambda}\alpha_3(\lambda_i).
    \end{equation}
    \item [Case $2$] $o(\lambda)=2t-1$, and $f_\lambda(1)\geq 4$. In this case, set $\alpha_3(\overline{2t-1})=\overline{2t+3}$. Define,
    \begin{equation}
        \mu=\zeta_3(\lambda)=\bigcup_{\lambda_i\in\lambda}\alpha_3(\lambda_i)\setminus(1^4).
    \end{equation}
    \item [Case $3$]  $o(\lambda)=2t-1$, $f_\lambda(1)\leq 3$ and $f_\lambda(2t+2)\geq 2$. In this case, set $\alpha_3(\overline{2t-1})=\overline{2t+3}$. Define
    \begin{equation}
        \mu=\zeta_3(\lambda)=\bigcup_{\lambda_i\in\lambda}\alpha_3(\lambda_i)\cup(t,t,t,t)\setminus(2t+2,2t+2).
    \end{equation}
    \item [Case $4$]  $o(\lambda)=2t-1$, $f_\lambda(1)=0$ and $f_\lambda(2t+2)=0$. In this case, set $\alpha_3(\overline{2t-1})=(t,t-3,\overline{2})$.  Define,
    \begin{equation}
    \mu=\zeta_3(\lambda)=\bigcup_{\lambda_i\in\lambda}\alpha_3(\lambda_i).
    \end{equation}
    \item[Case $5$]  $o(\lambda)=2t-1$, $f_\lambda(1)=1$ or $2$ and $f_\lambda(2t+2)=0$. In this case, set $\alpha_3(\overline{2t-1})=\overline{2t}$. Define,
    \begin{equation}
       \mu=\zeta_3(\lambda)=\bigcup_{\lambda_i\in\lambda}\alpha_3(\lambda_i)\setminus(1). 
    \end{equation}
    
    \item [Case $6$]  $o(\lambda)=2t-1$, $f_\lambda(1)=3$ and $f_\lambda(2t+2)=0$. In this case, set $\alpha_3(\overline{2t-1})=(t,t,\overline{2})$. Define,
    \begin{equation}
        \mu=\zeta_3(\lambda)=\bigcup_{\lambda_i\in\lambda}\alpha_3(\lambda_i)\setminus(1^3).
    \end{equation}

    \item[Case $7$]  $o(\lambda)=2t-1$, $f_\lambda(1)=0$ and $f_\lambda(2t+2)=1$. In this case, set $\alpha_3(\overline{2t-1})=(\overline{2t},t,t,1)$. Define,
    \begin{equation}
        \mu=\zeta_3(\lambda)=\bigcup_{\lambda_i\in\lambda}\alpha_3(\lambda_i)\setminus(2t+2).
    \end{equation}
    \item[Case $8$]  $o(\lambda)=2t-1$, $f_\lambda(1)=1$, $2$ or $3$ and $f_\lambda(2t+2)=1$. In this case, set $\alpha_3(\overline{2t-1})=(2t,2t,\overline{2})$. Define,
    \begin{equation}
    \mu=\zeta_3(\lambda)=\bigcup_{\lambda_i\in\lambda}\alpha_3(\lambda_i)\setminus(2t+2,1).
    \end{equation}
\end{itemize} 

We proceed to show that the image sets of the first three cases are in $E_t(n)$, and the last five cases are in $A_t^2(n)$. On the one hand, each $\mu$ in the first three cases satisfy that $o(\mu)=2t+1$ or $2t+3$. Moreover, each part is either $2t+2$ or cannot be divided by $t+1$. This yields $\mu\in E_t(n)$.  On the other hand, each $\mu$ in the last five cases satisfy that $o(\mu)=2,t$ or $2t$, which is even number and is not multiple of $t+1$. Moreover, from the construction of $\zeta_3$, it is easy to check that $f_\mu(2t+2)=0$, which implies $\mu$ is $t+1$-regular partition in all the five cases. Furthermore, it is clear that in each case, we have $o(\mu)=2t^k$ and $f_\mu(1)-\sum_{\mu_i\in\mu}g(\mu_i)$ in each case can be direct calculated, as listed in Table \ref{tb-2}. Thus we show that the image sets of the last five cases are in $A_t^2(n)$.
\begin{table}[]
    \centering
    \begin{tabular}{c|c}
    \hline
    & $f_{\mu}(1)- \sum_{\mu_i} g(\mu_i)$\\
    \hline
    Case 4 & $\leq-1$\\
    \hline
    Case 5 & $\leq-1$\\
    \hline
    Case 6 & $\leq-2$\\
    \hline
    Case 7 & $\leq-3$\\
    \hline
    Case 8 & $\leq-2$\\
    \hline
    \end{tabular}
    \caption{List of $f_{\mu}(1)- \sum_{\mu_i} g(\mu_i)$ in Case 4 $\sim$ Case 8.}
    \label{tb-2}
\end{table}

We next show that the image set of the $8$ cases are pairwise non-intersect. It is clear that the first three cases image are in $E_t(n)$ and the last five cases image are in $A^2_{t}(n)$. The details of these eight image sets are shown in Table \ref{tab:my_label1}.
\begin{table}[]
    \centering
    \begin{tabular}{c|c}
    \hline
    & $E_{t}(n)$\\
    \hline
    Case 1     &  $o(\mu)=2t+1$\\
    \hline
    Case 2 &  $o(\mu)=2t-1$, $ f_\mu(1)\geq \sum_{\mu_j\not=2t}g(\mu_j)$\\
    \hline
    Case 3 &       $o(\mu)=2t-1$, $ f_\mu(1)< \sum_{\mu_j\not=2t}g(\mu_j)$  \\
    \hline
    & $A_{t}^2(n)$\\
    \hline
    Case 4 & $o(\mu)=2$, $f_\mu(1)-\sum_{\mu_j\not=2t}g(\mu_j)=-1$\\
    \hline
    Case 5 & $o(\mu)=2t$, $f_\mu(1)-\sum_{\mu_j\not=2t}g(\mu_j)=0$ or $-1$\\
    \hline
    Case 6 & $o(\mu)=2$, $f_\mu(1)-\sum_{\mu_j\not=2t}g(\mu_j)=-2$\\
    \hline
    Case 7 & $o(\mu)=2t$, $f_\mu(1)-\sum_{\mu_j\not=2t}g(\mu_j)=-1$\\
    \hline
    Case 8 & $o(\mu)=2$, $f_\mu(1)-\sum_{\mu_j\not=2t}g(\mu_j)=0, 1$ or $2$\\
    \hline
    
    \end{tabular}
    \caption{The image sets of nine cases in Lemma \ref{zeta3}}
    \label{tab:my_label1}
\end{table}

We next show that $\zeta_3$ is an injection. To this end, let $I_{\zeta_3}(n)$ denote that the image set of $\zeta_3$, which is a subset of $A_{t}^2(n) \cup E_{t}(n)$. We proceed to establish a map $\eta_3$ from $I_{\zeta_3}(n)$ into $F_t(n)$, such that for any $\lambda\in F_t(n)$, 
\begin{equation}\label{eta_3}
    \eta_3(\zeta_3(\lambda))=\lambda.
\end{equation}
This implies $\zeta_3$ is an injection.

We now describe $\eta_3$. For any partition $\mu\in I_{\zeta_3}(n)$, we first define a map $\beta_3$. For the non-overlined part $\mu_j$, define
\begin{equation}
    \beta_3(\mu_j)=\begin{cases}
        \mu_j+g(\mu_j)&\text{if } \mu_j\ne 2t;\\
        \mu_j&\text{if }\mu_j=2t.
    \end{cases}
\end{equation}

For the overlined part $\overline{o(\mu)}$, if $\mu\in E_t(n)$, there are three cases.

\begin{itemize}
    \item [Case $1$] $o(\mu)=2t+1$, it is easy to check that $ f_\mu(1)\geq \sum_{\mu_j=\mu}g(\mu_j)$. Define $\beta_3(\overline{2t+1})=\overline{2t+1}$ and 
    \begin{equation}
        \lambda=\eta_3(\mu)=\bigcup_{\mu_i\in\mu}\beta_3(\mu_i)\setminus(1^{\sum_{\mu_j\not=2t}g(\mu_j)}).
    \end{equation}
    \item [Case $2$] $o(\mu)=2t+3$ and $ f_\mu(1)\geq \sum_{\mu_j\not=2t}g(\mu_j)$. Define $\beta_3(\overline{2t+3})=(\overline{2t-1},1^4)$ and 
    \begin{equation}
        \lambda=\eta_3(\mu)=\bigcup_{\mu_i\in\mu}\beta_3(\mu_i)\setminus(1^{\sum_{\mu_j\not=2t}g(\mu_j)}).
    \end{equation}
    \item [Case $3$] $o(\mu)=2t+3$ and $ f_\mu(1)< \sum_{\mu_j\not=2t}g(\mu_j)$. By Table \ref{tab:my_label1}, we see that this case is in Case $3$, which implies $f_\mu(t)\geq4$. Set $\nu=\mu\setminus(t^4)$. By the definition of $\beta_3$ and $\nu$, define $\beta_3(\overline{2t+3})=((2t+2)^2,\overline{2t-1})$ and 
    \begin{equation}
    \lambda=\eta_3(\mu)=\bigcup_{\nu_i\in\nu}\beta_3(\nu_i)\setminus(1^{\sum_{\mu_j\not=2t}g(\nu_j)}).
    \end{equation}
\end{itemize}

If $\mu\in A_t^2(n)$, then there are five cases for $\beta_3(\overline{o(\mu)})$.

\begin{itemize}
    \item [Case $1$] If $o(\mu)=2$ and $f_\mu(1)-\sum_{\mu_j\not= 2t}g(\mu_j)=-1$. Then by Table \ref{tab:my_label1}, we have $\mu$ is in the image set of Case $4$. By the construction of $\zeta_3$, we see that $f_\mu(t)\ge 1$ and $f_\mu(t-3)\ge 1$. Set $\nu=\mu\setminus(t,t-3)$. We also define $\beta_3(\overline{2})=\overline{2t-1}$ and
    \begin{equation}
        \lambda=\eta_3(\mu)=\bigcup_{\nu_i\in\nu}\beta_3(\nu_i)\setminus(1^{\sum_{\nu_j\ne 2t}g(\nu_j)}).
    \end{equation}
    \item [Case $2$] If $o(\mu)=2t$, and $f_\mu(1)-\sum_{\mu_j\not= 2t}g(\mu_j)=0$ or $1$. In this case, by Table \ref{tab:my_label1}, we have $\mu$ is in the image set of Case $5$. Set $\beta(\overline{2t})=(\overline{2t-1},1)$. Define, 
    \begin{equation}
        \lambda=\eta_3(\mu)=\bigcup_{\mu_i\in\mu}\beta_3(\mu_i)\setminus(1^{\sum_{\mu_j\ne 2t}g(\mu_j)}).
    \end{equation}
    \item [Case $3$] If $o(\mu)=2$ and $f_\mu(1)-\sum_{\mu_j\not= 2t}g(\mu_j)=-2$. In this case, by Table \ref{tab:my_label1}, we have $\mu$ is in the image set of Case $6$. And by the construction of $\zeta_3$, we see that $f_{\mu}(t)\ge 2$. Set $\nu=\mu\setminus(t,t)$ and $\beta_3(\overline{2})=(\overline{2t-1},1^3)$. Define, 
    \begin{equation}
        \lambda=\eta_3(\mu)=\bigcup_{\nu_i\in\nu}\beta_3(\nu_i)\setminus(1^{\sum_{\nu_j\ne 2t}g(\nu_j)}).
    \end{equation}
    \item [Case $4$] If $o(\mu)=2t$, $f_\mu(1)-\sum_{\mu_j\not= 2t}g(\mu_j)=-1$. By the construction of $\zeta_3$, we see that $f_\mu(t)\geq2$ and $f_\mu(1)\geq1$. Set $\nu=\mu\setminus(t,t,1)$ and $\beta_3(\overline{2t})=(\overline{2t-1},2t+2)$. Define,
    \begin{equation}
        \lambda=\eta_3(\mu)=\bigcup_{\nu_i\in\nu}\beta_3(\nu_i)\setminus(1^{\sum_{\nu_j\ne 2t}g(\nu_j)}).
    \end{equation}
    \item[Case $5$] If $o(\mu)=2$, $f_\mu(1)-\sum_{\mu_j\not= 2t}g(\mu_j)\geq0$. By the construction of $\zeta_3$, we see that $f_\mu(2t)\geq2$. Set $\nu=\mu\setminus(2t,2t)$ and $\beta_3(\overline{2})=(\overline{2t-1},2t+2,1)$. Define,
    \begin{equation}
         \lambda=\eta_3(\mu)=\bigcup_{\nu_i\in\nu}\beta_3(\nu_i)\setminus(1^{\sum_{\nu_j\ne 2t}g(\nu_j)}).
    \end{equation}
\end{itemize}

It is routine to check $\eta_3$ satisfies \eqref{eta_3}. Thus $\zeta_3$ is an injection when $t\ge 6$. 

We now consider the case $t=4$. Actually, the proof of Lemma \ref{zeta3} are still valid for the case $t=4$ except for the Case $4$. The invalid case is because that $\alpha_3(\overline{7})=(4,\overline{2},1)$, which increase the number of $f_\mu(1)$. So we need to find a new image set in $\hat{A}_8^4(n)$. Now we modify this case as follows:
\begin{equation}\label{eq-def-alpha3-4}
    \alpha_3(\overline{7})=\begin{cases}
   (\overline{8}')& \text{if there exists }\lambda_i \text{ such that }5\mid \lambda_i;\\
    (\overline{4}',{3})& \text{otherwise.}
\end{cases}
\end{equation}
Moreover, the map $\zeta_3$ also need to modified as follows:
\begin{equation}\label{eq-def-alpha3-4-1}
    \mu=\zeta_3(\lambda)=\begin{cases}
    \bigcup_{\lambda_i\in\lambda}\alpha_3(\lambda_i)\setminus(1) & \text{if there exists }\lambda_i \text{ such that }5\mid \lambda_i;\\
    \bigcup_{\lambda_i\in\lambda}\alpha_3(\lambda_i) & \text{otherwise.}
    \end{cases}
\end{equation}

If there exists $\lambda_i$ such that $5\mid \lambda_i$, we see that $o(\mu)=8=2\cdot4$. Moreover, it is clear that $f_\mu(1)-\sum_{\mu_i\not=8}g(\mu_i)=-1$ which implies that $f_\mu(1)-\sum_{\mu_i}g(\mu_i)<0$. This implies that $\mu\in \hat{A}_4^4(n)$. Otherwise, we can find that $o(\mu)=4=1\cdot4$. Moreover, it is clear that $f_\mu(1)=0<g(o(\mu))=1$. This implies that $\mu\in \hat{A}_4^4(n)$.

We next show that $\zeta_3$ is an injection in Case $4$. It suffices to show that the restrictive map of $\zeta_3$ in Case 4 to $\hat{A}_4^4(n)$ is an injection. To this end, let $I_4(n)=I_{\zeta_3}(n)\cap \hat{A}_4^4(n)$ denote the image set of $\zeta_3$ restrict to Case 4. We proceed to establish a map $\eta_3$ from $I_4(n)$ into $F_t(n)$, such that for any $\lambda$ in Case 4, we have
\begin{equation}\label{equ-eta3-case 4}
    \eta_3(\zeta_3(\lambda))=\lambda.
\end{equation}
This implies $\zeta_3$ is an injection in Case 4.

We now describe $\eta_3$ in Case 4. For any partition $\mu\in I_4(n)$, we first define a map $\beta_3$. For the non-overlined part $\mu_j$, define
\begin{equation}
    \beta_3(\mu_j)=\begin{cases}
        \mu_j+g(\mu_j)&\text{if } \mu_j\ne 8;\\
        \mu_j&\text{if }\mu_j=8.
    \end{cases}
\end{equation}

For the overliend part $\overline{o(\mu)}'$, from the construction of \eqref{eq-def-alpha3-4}, we see $o(\mu)=8$ or $4$. There are two cases.

\begin{itemize}
    \item [Case 1] If $o(\mu)=8$, set $\beta_3(\overline{8}')=(\overline{7},1)$, by definition \eqref{eq-def-alpha3-4-1}, we see that $f_\mu(1)=\sum_{\mu_j\ne 8}g(\mu_j)-1$. Hence we may define
    \begin{equation}
        \lambda=\eta_3(\mu)=\bigcup_{\mu_i\in\mu}\beta_3(\mu_i)\setminus(1^{\sum_{\mu_j\ne 8}g(\mu_j)}).
    \end{equation}
    \item [Case 2] If $o(\mu)={4}$, by definition \eqref{eq-def-alpha3-4}, we have $5\nmid \mu_i$ for each $\mu_i$. Moreover, $f_\mu(3)\ge 1$. Set $\nu=\mu\setminus(3)$ and $\beta_3(\overline{4}')=\overline{7}$, define
    \begin{equation}
        \lambda=\eta_3(\mu)=\bigcup_{\nu_i\in\nu}\beta_3(\nu_i).
    \end{equation}
\end{itemize}

It is routine to check $\eta_3$ satisfies \eqref{equ-eta3-case 4}. Thus $\zeta_3$ is an injection when $t=4$.
\qed

For example, set $t=6$ and $\lambda^1=(14,\overline{13},12,4,2,1)$. Applying the injection $\zeta_3$, it is easy to check that $\mu^1=\zeta_3(\lambda^1)=(\overline{13},12,12,4,2,1^3)$. For another example, set $t=6$ and $\lambda^2=(\overline{11},4,2,1^5)$, it is easy to check that $\mu^2=\zeta_3(\lambda^2)=(\overline{15},4,2,1)$; $\lambda^3=(\overline{11},5,2)$, it is easy to check that $\mu^3=\zeta_3(\lambda^3)=(6,5,3,\overline{2},2)$. Moreover, it is routine to check that $\eta_3(\mu)=\lambda$, when $\mu=\mu^1$, $\mu^2$ or $\mu^3$.

We are now in a position to show Theorem \ref{thm-com-int-umn} for even number $t$.


{\noindent\it Proof of Theorem \ref{thm-com-int-umn} for even number $t$.} From Lemma \ref{zeta1} $\sim$ \ref{zeta3}, we see that \eqref{bt+12-bt2} has non-negative coefficients when $t\geq4$ even. Thus it proves Theorem \ref{thm-com-int-umn} for even number $t$. \qed

\noindent{\bf Acknowledgments.}   This work was supported by the National Science Foundation of China grants 12171358 and 12371336.

\end{document}